\documentclass{amsart}

\usepackage{amsthm,amsfonts,amsmath,amssymb,latexsym,epsfig,mathrsfs,yfonts,marvosym}
\usepackage[usenames]{color}

\DeclareMathAlphabet\oldmathcal{OMS}        {cmsy}{b}{n}
\SetMathAlphabet    \oldmathcal{normal}{OMS}{cmsy}{m}{n}
\DeclareMathAlphabet\oldmathbcal{OMS}       {cmsy}{b}{n}
 
\usepackage{eucal}

\usepackage[active]{srcltx}
\usepackage{graphicx}
\usepackage[all]{xy}
\usepackage{epsfig}
\usepackage{layout}

\newtheorem{theorem}{Theorem}[section]
\newtheorem{lemma}[theorem]{Lemma}
\newtheorem{proposition}[theorem]{Proposition}

\newtheorem{definition}[theorem]{Definition}
\newtheorem{definition/proposition}[theorem]{Definition/Proposition}

\newenvironment{remark}{\medskip \refstepcounter{theorem}
\noindent  {\bf Remark \thetheorem}.\rm}{\,}

\def\BOne{{\mathchoice {\rm 1\mskip-4mu l} {\rm 1\mskip-4mu l}
                          {\rm 1\mskip-4.5mu l} {\rm 1\mskip-5mu l}}}

\def\fract#1#2{\raise4pt\hbox{$ #1 \atop #2 $}}

\def\bbc{{\mathbb C}}

\def\bbp{{\mathbb P}}

\def\bbr{{\mathbb R}}

\def\bbz{{\mathbb Z}}

\def\gra{\alpha}

\def\grg{\gamma}

\def\grk{\kappa}

\def\gro{\omega}

\def\grz{\zeta}

\def\cala{{\mathcal A}}
\def\calb{{\mathcal B}}
\def\calc{{\mathcal C}}

\def\cald{{\mathcal D}}

\def\calf{{\mathcal F}}

\def\cals{{\oldmathcal S}}

\def\calw{{\mathcal W}}

\def\calz{{\oldmathcal Z}}

\def\la#1{\hbox to #1pc{\leftarrowfill}}
\def\ra#1{\hbox to #1pc{\rightarrowfill}}

\def\ge{{\mathfrak e}}

\def\gm{{\mathfrak m}}
\def\gn{{\mathfrak n}}
\def\go{{\mathfrak o}}

\def\gt{{\mathfrak t}}

\def\gy{{\mathfrak y}}

\def\gB{{\mathfrak B}}
\def\gC{{\mathfrak C}}

\def\gE{{\mathfrak E}}

\def\gR{{\mathfrak R}}
\def\gS{{\mathfrak S}}

\def\teta{\tilde{\eta}}

\def\teta{\tilde{{\eta}}}

\def\tgro{\tilde{\gro}}
\def\tcald{\tilde{\cald}}
\def\teta{\tilde{\eta}}

\def\cJ{\hat{J}}

\def\hook{\mathbin{\hbox to 6pt{%
                 \vrule height0.4pt width5pt depth0pt
                 \kern-.4pt
                 \vrule height6pt width0.4pt depth0pt\hss}}}

\begin{document}
\bibliographystyle{amsalpha}

\title{Extremal Sasakian Metrics on ${\bf S^3}$-bundles over ${\bf S^2}$}\thanks{In the beginning of this work the author was partially supported by NSF grant DMS-0504367.}

\author{Charles P. Boyer}
\address{Department of Mathematics and Statistics,
University of New Mexico, Albuquerque, NM 87131.}

\email{cboyer@math.unm.edu} 

\begin{abstract}
In this note I prove the existence of extremal Sasakian metrics on $S^3$-bundles over $S^2.$ These occur in a collection of open cones that I call a bouquet.
\end{abstract}

\maketitle

\section{Introduction}
It is well known that the $S^3$-bundles over $S^2$ are classified by $\pi_1(SO(4))=\bbz_2$. So there are exactly two such bundles, the trivial bundle $S^2\times S^3$ and one non-trivial bundle  $X_\infty$ (in Barden's notation \cite{Bar65}). They are distinquished by their second Stiefel-Whitney class $w_2\in H^2(M,\bbz_2)$. 

There are infinitely many Sasaki-Einstein metrics\footnote{For basic material concerning Sasakian geometry I refer to our recent book \cite{BG05}.} on $S^2\times S^3$ which belong to a toric contact structure \cite{GMSW04a,GMSW04b,CLPP05,MaSp05b}. These, of course, are all extremal, but they have $c_1(\cald)=0$ where $\cald$ is the contact bundle. There is a more or less obvious constant scalar curvature extremal Sasakian metric on $S^2\times S^3$ which is not Sasaki-Einstein and has first Chern class $c_1(\cald)=2(k_1-k_2)\gra$ for every pair $(k_1,k_2)$ of relatively prime positive integers. Here $\gra$ is a generator of $H^2(S^2\times S^3,\bbz)$, and without loss of generality we can assume that $k_1> k_2$. We recover the well known Kobayashi-Tanno homogenous Einstein metric in the case $(k_1,k_2)=(1,1)$. The Sasakian structures are constructed from the K\"ahler form $\gro_{k_1,k_2}=k_1\gro_1+k_2\gro_2$ on $S^2\times S^2$ with the product complex structure, where $\gro_1(\gro_2)$ are the standard symplectic forms on the first (second) factor. The metric corresponding to this K\"ahler form has constant scalar curvature. One then forms the $S^1$-bundle over $S^2\times S^2$ whose cohomology class is $[\gro_{k_1,k_2}]$. The constant scalar curvature K\"ahler metrics lifts to a constant scalar curvature Sasakian metric on $S^2\times S^3$ which is homogeneous, hence toric. So there is a three dimensional Sasaki cone $\grk$ as described in \cite{BGS07b}, and by the openness theorem of \cite{BGS06} there is an open set of extremal Sasakian metrics in $\grk$. I want to emphasize that although extremal quasi-regular Sasakian metrics are always lifts of extremal K\"ahlerian orbifold metrics, it is NOT true that the openness theorem for extremal Sasaki metrics is obtained by simply lifting the openness theorem \cite{LeSi93b,LeSi94} for extremal K\"ahler metrics. In this note I shall prove that there are many other extremal Sasakian metrics on $S^2\times S^3$ belonging to the same contact structure.

In contrast to the situation of $S^2\times S^3$, until now there are no known extremal Sasakian metrics on $X_\infty$. Here I also show that $X_\infty$ admits many extremal Sasakian structures belonging to the same contact structure.

\section{Review of Extremal Sasakian Metrics}
Recall \cite{BG05} that a Sasakian structure $\cals=(\xi,\eta,\Phi,g)$ on a smooth manifold $M$ is a contact metric structure such that $\xi$ preserves the underlying almost CR structure $(\cald,J)$ defined by $\cald=\ker \eta$ and $J=\Phi|_\cald$, and the almost CR structure is integrable. Now we deform the contact 1-form by $\eta\mapsto \eta(t)=\eta+t\grz$ where $\grz$ is a basic 1-form with respect to the characteristic foliation $\calf_\xi$ defined by the Reeb vector field $\xi.$ Here $t$ lies in a suitable interval containing $0$ and such that $\eta(t)\wedge d\eta(t)\neq 0$. This gives rise to a family of Sasakian structures $\cals(t)=(\xi,\eta(t),\Phi(t),g(t))$ that we denote by ${\mathcal S}(\xi, \bar{J})$ where $\bar{J}$ is the induced complex structure on the quotient bundle $\nu(\calf_\xi)=TM/L_\xi$ by the trivial line bundle generated by $\xi.$ As the notation suggests we always assume that $\cals(0)=(\xi,\eta(0),\Phi(0),g(0))=\cals$.

Note that this deforms the contact structure (hence, the CR structure) from $\cald$ to $\cald_t=\ker \eta(t)$, but they are isomorphic as complex vector bundles and isotopic as contact structures by Gray's theorem (cf. \cite{BG05}). In fact, each choice of 1-form defines a splitting 
$TM=L_\xi+\cald$ and an isomorphism $\cald\approx \nu(\calf_\xi)$ of complex vector bundles. The complex structure $J$ on $\cald$ defines a further splitting of the complexified bundle $\cald\otimes \bbc=\cald^{1,0}+\cald^{0,1}$, and the usual Dolbeault type complexes with transverse Hodge theory holds \cite{BG05}, and the same for $\cald_t$.

We assume that $M$ is compact of dimension $2n+1$ with a Sasakian structure $\cals$, and note that the associated Riemannian metric is uniquely determined by $\eta$ and $\Phi$ as
$$g=d\eta\circ (\Phi\otimes \BOne) \oplus \eta\otimes \eta.$$ 
Following \cite{BGS06} we let $s_g$ denote the scalar curvature of $g$ and define 
the ``energy functional'' $E:{\mathcal S}(\xi,\bar{J})\ra{1.4} \bbr$ by
\begin{equation}\label{var}
E(g) ={\displaystyle \int _M s_g ^2 d{\mu}_g ,}\, 
\end{equation}
i.e. the $L^2$-norm of the scalar curvature. Critical points $g$ of this functional are called {\it extremal Sasakian metrics}. In this case we also say that the Sasakian structure $\cals$ is extremal. Similar to the K\"ahlerian case, the Euler-Lagrange equations for this functional give \cite{BGS06}
\begin{theorem}\label{ELeqn}
A Sasakian structure $\cals\in {\mathcal S}(\xi,\bar{J})$ is a critical point for the energy functional (\ref{var}) if and only if the gradient vector field $\partial^\#_gs_g$ is transversely holomorphic. In particular, Sasakian metrics with constant scalar curvature are extremal.
\end{theorem}

Here $\partial^\#_g$ is the $(1,0)$-gradient vector field defined by $g(\partial^\#_g\varphi,\cdot)= \bar{\partial}\varphi$.
It is important to note that a Sasakian metric $g$ is extremal if and only if the `transverse metric' $g^T=d\eta\circ (\Phi\otimes \BOne)$ is extremal in the K\"ahlerian sense which follows from the well known relation between scalar curvatures $s_g=s_g^T-2n$ where $s_g^T$ is the scalar curvature of the transverse metric. It follows that

\begin{proposition}\label{SasKah}
Let $\cals=(\xi,\eta,\Phi,g)$ be a Sasakian structure on $M$ of dimension $2n+1$ and let $U\subset M$ be an open set such that $f:U\ra{1.5} \bbc^n$ is a local submersion. Then the restriction $g|_U$ is an extremal Sasakian metric if and only if $g^T$ viewed as a K\"ahlerian metric on $f(U)$ is extremal. In particular, if the Sasakian structure $\cals$ is quasi-regular and $(\gro,\cJ,h)$ is the induced orbifold K\"ahlerian structure on the quotient $\calz=M/\calf_\xi$, then $g$ is Sasakian extremal if and only if $h$ is K\"ahlerian extremal. Moreover, $g$ has constant scalar curvature if and only if $h$ has constant scalar curvature.
\end{proposition}

Given a strictly pseudoconvex CR structure $(\cald,J)$ of Sasaki type on a smooth manifold $M$ of dimenision $2n+1$, we consider the set ${\mathcal S}(\cald,J)$ of Sasakian structures whose underlying CR structure is $(\cald,J)$. The group $\gC\gR(\cald,J)$ of CR transformation acts on ${\mathcal S}(\cald,J)$, and the quotient space $\grk(\cald,J)$ is called the Sasaki cone \cite{BGS06} whose dimension satisfies $1\leq \dim \grk(\cald,J)\leq n+1$.  For a strictly pseudoconvex CR structure $(\cald,J)$ on a compact manifold the group $\gC\gR(\cald,J)$ is compact except in the case of standard CR structure on the sphere $S^{2n+1}$ \cite{Sch95,Lee96} in which case it is isomorphic to $SU(n+1,1)$ \cite{Web77}. Thus, $\gC\gR(\cald,J)$ has a unique maximal torus up to conjugacy.  So, as discussed in \cite{Boy10a} it is often convenient to consider the `unreduced' Sasaki cone $\gt^+=\gt^+(\cald,J)$ where $\gt$ is the Lie algebra of a maximal torus in the group $\gC\gR(\cald,J)$ of CR transformations. This is defined by
\begin{equation}\label{uSascone}
\gt^+=\{\xi'\in \gt~|~\eta(\xi')>0\}
\end{equation}
where $\eta$ is any 1-form representing $\cald$, and is an open convex cone in $\gt$. It is related to the Sasaki cone $\grk(\cald,J)$ by $\grk(\cald,J)=\gt^+(\cald,J)/\calw$ where $\calw$ is the Weyl group of $\gC\gR(\cald,J)$. By abuse of terminology I also refer to $\gt^+$ as the Sasaki cone. Note that with $(\cald,J)$ fixed, choosing a Reeb field uniquely chooses a 1-form $\eta$ such that $\ker\eta=\cald$, and with $J$ also fixed $\Phi$, hence, $g$ are uniquely specified. Thus, we can think of the Sasaki cone $\gt^+(\cald,J)$ as consisting of Sasakian structures $\cals=(\xi,\eta,\Phi,g)$, so again by abuse of notation we can write $\cals\in \gt^+(\cald,J)$.

Let $\eta(t)=\eta+t\grz$ where $\grz$ is a basic 1-form that is invariant under the full torus $T$. Such a 1-form can always be obtained by averaging over $T$. So the Lie algebra $\gt$ does not change under such a deformation, but generally the Sasaki cone $\gt^+$ associated with $\cald_t$ shifts, and we denote it by $\gt^+(t)$. We let ${\mathcal S}^T(\xi,\bar{J})$ denote the subset of ${\mathcal S}(\xi,\bar{J})$ consisting of $T$-invariant Sasakian structures. If $\cals=(\xi,\eta,\Phi,g)$ is $T$-invariant, then ${\mathcal S}^T(\xi,\bar{J})$ consists of all deformations obtained by $\eta\mapsto \eta(t)=\eta+t\grz$ with $\grz$ invariant under $T$. We are interested in the case when the torus $T$ has maximal dimension and the Reeb vector field is an element of the Lie algebra $\gt$, that is, a toric contact structure $\cald$ of {\it Reeb type}. 

\begin{lemma}\label{invlem}
Let $(M,\cald,T)$ be a toric contact structure of Reeb type with Reeb vector field $\xi.$ Suppose also that there is an extremal representative $\cals(t)=(\xi,\eta(t),\Phi_t,g_t)\in {\mathcal S}(\xi,\bar{J})$. Then $\cals(t)\in {\mathcal S}^T(\xi,\bar{J})$.
\end{lemma}

\begin{proof}
Since the contact structure $\cald$ is toric of Reeb type, there is a compatible $T$-invariant Sasakian structure $\cals=(\xi,\eta,\Phi,g)$ by \cite{BG00b}, and suppose that ${\mathcal S}(\xi,\bar{J})$ has an extremal representative $\cals(t)=(\xi,\eta(t),\Phi_t,g_t)$. By the Rukimbira Approximation Theorem, if necessary, there is a quasi-regular $T$-invariant Sasakian structure invariant close to $\cals$, so we can assume that $\cals$ is quasi-regular. But then by Proposition \ref{SasKah} the Sasakian deformation corresponds to a K\"ahlerian deformation on the base K\"ahler orbifold $M/\calf_\xi$. By a theorem of Calabi \cite{Cal85} this deformed K\"ahler structure has maximal symmetry, and one easily sees that the corresponding deformed Sasakian structure $\cals(t)$ also has maximal symmetry which implies that $\cals(t)\in {\mathcal S}^T(\xi,\bar{J})$.
\end{proof}

Now let $(\cald,J)$ be a strictly pseudoconvex CR structure of Sasaki type. We say that $\cals\in \gt^+(\cald,J)$ is an {\it extremal element}\footnote{In \cite{BGS06} this was called a canonical element, but I prefer to call it an extremal element.} of $\gt^+(\cald,J)$ if there exists an extremal Sasakian structure in ${\mathcal S}(\xi,\bar{J})$. We also say that $\cals$ is an extremal element of $\grk(\cald,J)$. So we can define the {\it extremal set} $\ge(\cald,J)\subset \grk(\cald,J)$ as the subset consisting of those elements that have extremal representatives (similarly for $\gt^+(\cald,J)$).

Recall the {\it transverse homothety} (cf. \cite{BG05}) taking a Sasakian structure $\cals=(\xi,\eta,\Phi,g)$ to the Sasakian structure $$\cals_a=(a^{-1}\xi,a\eta,\Phi,ag+(a^2-a)\eta\otimes \eta)$$
for any $a\in \bbr^+$. It is easy to see that

\begin{lemma}\label{transhomoext}
Let $\cals$ be a Sasakian structure. If $\cals$ is extremal, so is $\cals_a$, and if $\cals$ has constant scalar curvature so does $\cals_a$. In particular, if the extremal set $\ge(\cald,J)$ is non-empty it contains an extremal ray of Sasakian structures.
\end{lemma}

A main result of \cite{BGS06} says that $\ge(\cald,J)$ is an open subset of $\grk(\cald,J)$. Lemma \ref{transhomoext} says that $\ge(\cald,J)$ is conical in the sense that it is a union of open cones, so it is not necessarily connected.

\section{Bouquets of Sasaki Cones and Extremal Bouquets}

There may be many Sasaki cones associated to a given contact structure $\cald$ of Sasaki type. They are distinguished by their complex structures $J.$ As shown in \cite{Boy10a} to a compatible almost complex structure $J$ on a compact contact manifold $(M,\cald)$ one can associate a conjugacy class $\calc_T(\cald)$ of maximal tori in the contactomorphism group $\gC\go\gn(M,\cald)$. Furthermore, almost complex structures that are equivalent under a contactomorphism give the same conjugacy class $\calc_T(\cald)$. So given inequivalent complex structures $J_l$ labelled by positive integers, we can associate unreduced Sasaki cones $\gt^+(\cald,J_l)$, or the full Sasaki cones $\grk(\cald,J_l)$. This leads to
\begin{definition}\label{Sasbou}
We define a {\bf bouquet of Sasaki cones} $\gB(\cald)$ as the union 
$$\gB(\cald)=\cup_{l\in \cala}\grk(\cald,J_l)$$ 
where $\cala\subset \bbz^+$ is an ordered subset. We say the it is an {\bf $N$-bouquet} if the cardinality of $\cala$ is $N$ and denote it by $\gB_N(\cald)$.
\end{definition}

Clearly a $1$-bouquet of Sasaki cones is just a Sasaki cone. Generally, the Sasaki cones in $\gB(\cald)$ can have varying dimension.

The example of interest to us here has in its foundations in the work of Karshon \cite{Kar03} and Lerman \cite{Ler03b}. The general formulation is given in \cite{Boy10a}. We begin with a simply connected symplectic orbifold $(\calb,\gro)$ such that $\gro$ defines an integral class in the orbifold cohomology group $H^2_{orb}(\calb,\bbz)$. Suppose further that the class $[\gro]$ is primitive and that there are compatible almost complex structures $\cJ_l$ on $\calb$ such that $(\gro,\cJ_l)$ is K\"ahler for each $l$ in some index set $\cala.$ We can associate to each such $\cJ_l$ a conjugacy class of maximal tori in the symplectomorphism group $\gS\gy\gm(\calb,\gro)$. Assume that this map is injective. Now form the principal $S^1$-orbibundle over $\calb$ associated to $[\gro]$ and assume that the total space $M$ is smooth. By the orbifold version of the Boothby-Wang construction we get a contact manifold $(M,\eta)$ satisfying $\pi^*\gro=d\eta$ where $\pi:M\ra{1.5} \calb$ is the natural orbifold projection map. We can lift a compatible almost complex structure $\cJ_l$ on $\calb$ to a $\cald$-compatible almost complex structure $J_l$ on $\cald$. Choosing tori $T_l$ associated to $\cJ_l$ we can also lift these to maximal tori $\Xi\times \pi^{-1}(T_l)$ in the contactomorphism group $\gC\go\gn(M,\cald)$ where $\cald=\ker\eta$ and $\Xi$ is the circle group generated by the Reeb vector field $\xi$ of $\eta$. Suppose that there are exactly $N$ such maximal tori in $\gC\go\gn(M,\cald)$, then we get an N-bouquet $\gB_N(\cald)$ of Sasaki cones on $(M,\cald)$ which intersect in the ray of Reeb vector fields, $\xi_a=a^{-1}\xi$. If we projectivize the Sasaki cones by transverse homotheties we obtain the usual notion of bouquet, namely, a wedge product $\bigvee\bbp(\grk(\cald,\cJ_l))$ with base point $\xi$.

Suppose that $\cals_l=(\xi,\eta,\Phi_l,g_l)$ is an extremal element of $\grk(\cald,J_l)$ for each $l\in \cala$, that is for each $l\in \cala$ there is an extremal representative $\cals_l(t)\in {\mathcal S}(\xi,\bar{J}_l)$. Then by the openness theorem there is a nonempty open extremal set $\ge(\cald,J_l)\subset \grk(\cald,J_l)$ containing $\cals=(\xi,\eta,\Phi,g)$ for each $l\in \cala$. This leads to

\begin{definition}\label{extbouq}
We define an {\bf extremal bouquet} $\gE\gB(\cald)$ associated with the contact structure $\cald$ to be the union $\cup_{l\in \cala}\ge(\cald,J_l)$ if for each $l\in \cala$ the extremal set $\ge(\cald,J_l)$ is non-empty. Moreover, it is called an {\bf extremal $N$-bouquet} $\gE\gB_N(\cald)$ if $\cala$ has cardinality $N$.
\end{definition}

One easily sees from the openness theorem of \cite{BGS06} that $\gE\gB(\cald)$ is open in $\gB(\cald)$. An important open problem here is to obtain a good measure of the size of the extremal sets $\ge(\cald,J_l)$, and hence, a measure of the size of the extremal bouquet. The actual size of extremal Sasakian sets is known in very few cases, namely the standard sphere, and the Heisenberg group \cite{BGS06,Boy09}.

\section{The Main Theorems}

Let us now construct contact structures of Sasaki type on $S^2\times S^3$ and $X_\infty$. First consider the toric symplectic manifold $S^2\times S^2$ with symplectic form $\gro_{k_1,k_2}=k_1\gro_1+k_2\gro_2$ on $S^2\times S^2$ where $(k_1,k_2)$ are relatively prime integers satisfying $k_1\geq k_2$, and $\gro_1(\gro_2)$ is the standard symplectic forms on the first (second) factor of $S^2$, respectively. Let $\pi:M\ra{1.5} S^2\times S^2$ be the circle bundle corresponding to the cohomology class $[\gro_{k_1,k_2}]\in H^2(M,\bbz)$. As stated in the introduction $M$ is diffeomorphic to $S^2\times S^3$ for each such pair $(k_1,k_2)$. By the Boothby-Wang construction one obtains a contact structure $\cald_{k_1,k_2}$ on $M$ by choosing a connection 1-form $\eta_{k_1,k_2}$ such that $\pi^*\gro_{k_1,k_2}=d\eta_{k_1,k_2}$. Justin Pati \cite{Pat09} has recently shown using contact homology that there are infinitely many such contact structures with the same first Chern class that are inequivalent.

\begin{theorem}\label{thmA}
Let $M=S^2\times S^3$ with the contact structure $\cald_{k_1,k_2}$ described above. Let $N=\lceil \frac{k_1}{k_2} \rceil$ denote the smallest integer greater than or equal to $\frac{k_1}{k_2}$. Then for each pair $(k_1,k_2)$ of relatively prime integers satisfying $k_1\geq k_2$, there is an extremal $N$-bouquet $\gE\gB_N(\cald_{k_1,k_2})$ of toric Sasakian structures associated to the contact structure $\cald_{k_1,k_2}$ on $S^2\times S^3$, and all the cones of the bouquet have dimension three. Furthermore, the extremal Sasakian structure corresponding to the contact 1-form $\eta_{k_1,k_2}$ has constant scalar curvature if and only if the transverse complex structure is that induced by the product complex structure on the base $S^2\times S^2$.
\end{theorem}

\begin{proof}
By Proposition \ref{SasKah} a Sasakian structure $\cals=(\xi,\eta,\Phi,g)$ in $\grk(\cald_{k_1,k_2},J)$ is extremal if and only if the induced K\"ahler structure $(\gro_{k_1,k_2},\cJ,h)$ on $S^2\times S^2$ is extremal. Now Karshon \cite{Kar03} proved that there are precisely $\lceil \frac{k_1}{k_2} \rceil$ even Hirzebruch surfaces  $S_{2m}$ that are K\"ahler with respect to the symplectic form $\gro_{k_1,k_2}$. More precisely for each $m=0,\cdots, \lceil \frac{k_1}{k_2} \rceil-1$ there is a diffeomorphism that takes the K\"ahler form $S_{2m}$ to $\gro_{k_1,k_2}$. Note that $m=0$ corresponds the product complex structure $\bbc\bbp^1\times \bbc\bbp^1$. Moreover, Calabi \cite{Cal82} has shown that for every K\"ahler class $[\gro]$ of any Hirzebruch surface $S_{n}$ there is an extremal K\"ahler metric which is of constant scalar curvature if and only if $n=0$. Furthermore, as shown by Lerman \cite{Ler03b} one can lift maximal tori of the symplectomorphism group $\gS\gy\gm(S^2\times S^2,\gro_{k_1,k_2})$ to maximal tori in the contactomorphism group $\gC\go\gn(S^2\times S^3,\cald_{k_1,k_2})$ (see also Theorem 6.4 of \cite{Boy10a} for the more general situation). As explained above this gives a Sasaki cone associated to each of the transverse complex structures induced by the Hirzebruch surfaces, and hence an $N$-bouquet of Sasaki cones where $N= \lceil \frac{k_1}{k_2} \rceil$. Furthermore, since each Sasaki cone has an extremal Sasakian structure, namely the one with 1-form $\eta_{k_1,k_2}$, we can apply Theorem 7.6 of \cite{BGS06} to obtain an open set of extremal Sasakian structures in each of the Sasaki cones. This gives an extremal $N$-bouquet $\gE\gB_N(\cald_{k_1,k_2})$ as defined above.
\end{proof}

There is another family of contact structures on $S^2\times S^3$ that arise as circle bundles over the non-trivial $S^2$-bundle over $S^2$ which is diffeomorphic to $\bbc\bbp^2$ blown up at a point. Following \cite{Kar03} I denote this manifold by $\widetilde{\bbc\bbp}^2$. The construction of this family also gives contact structures on the non-trivial bundle $X_\infty$. Let $E$ be the exceptional divisor and let $L$ be a projective line in $\widetilde{\bbc\bbp}^2$ that does not intersect $E$. We consider the symplectic form $\tgro_{l,e}$ such that the symplectic areas of $L$ and $E$ are $2\pi l$ and $2\pi e$, respectively, where $(l,e)$ are relatively prime integers satisfying $l>e\geq 1$. Karshon \cite{Kar03} proved that there are precisely $\lceil \frac{e}{l-e} \rceil$ odd Hirzebruch surfaces $S_{2m+1}$ with $m=0,\cdots, \lceil \frac{e}{l-e} \rceil-1$ that are K\"ahler with respect to the symplectic form $\tgro_{l,e}$, and again there are diffeomorphisms $\psi_m$ that take the K\"ahler forms of $S_{2m+1}$ to the form $\tgro_{l,e}$. As occured in the proof of the previous theorem the different Hirzebruch surfaces correspond to non-conjugate maximal tori in the symplectomorphism group $\gS\gy\gm(\widetilde{\bbc\bbp}^2,\tgro_{l,e})$. By Proposition 6.3 and Theorem 6.4 of \cite{Boy10a} these can be lifted to non-conjugate tori in the contactmorphism group $\gC\go\gn(M,\tcald_{l,e})$ where $M$ is the circle bundle over $\widetilde{\bbc\bbp}^2$ corresponding to the cohomology class $[\gro_{l,e}]$ and $\tcald_{l,e}$ is the contact structure induced by a connection 1-form $\eta_{l,e}$ satisfying $\pi^*\gro_{l,e}=d\eta_{l,e}$. I assume that $l$ and $e$ are relatively prime positive integers. Now in order to identify the diffeomorphism type of $M$ it is enough by the Barden-Smale classification of simply connected 5-manifolds to compute the mod 2 reduction of the first Chern class $c_1(\tcald_{l,e})$. For this we pullback $c_1(\widetilde{\bbc\bbp}^2)$ to $M$ and use the relation $[\pi^*\tgro_{l,e}]=0$. Now $c_1(\widetilde{\bbc\bbp}^2)=2\gra_E+\gra_L$ where $\gra_E(\gra_L)$ is Poincar\'e dual of the divisor class $E(L)$, respectively  \cite{GrHa78}. We find $c_1(\tcald_{l,e})=\bigl(l-2e\bigr)\grg$ where $\grg$ is a generator of $H^2(M,\bbz)\approx \bbz$. It follows that $M$ is $S^2\times S^3$ if $l$ is even and $X_\infty$ if $l$ is odd. Then, using Calabi \cite{Cal82}, just as in the proof of Theorem \ref{thmA} we arrive at

\begin{theorem}\label{thmB}
Let $M$ be the circle bundle over $\widetilde{\bbc\bbp}^2$ with the contact structure $\tcald_{l,e}$ as described above, and let $N=\lceil \frac{e}{l-e} \rceil$. Then for each pair $(l,e)$ of relatively prime integers satisfying $l>e\geq 1$, there is an extremal $N$-bouquet $\gE\gB_N(\tcald_{l,e})$ of toric Sasakian structures associated to the contact structure $\tcald_{l,e}$, and the cones of the bouquet all have dimension three.  Furthermore, if $l$ is even $M=S^2\times S^3$, whereas, if $l$ is odd $M=X_\infty$ and in either case the extremal Sasakian structure corresponding to the contact 1-form $\eta_{l,e}$ does not have constant scalar curvature. 
\end{theorem}

\begin{remark}
If the pairs of integers $(k_1,k_2)$ and $(l,e)$ are not relatively prime, but $(k_1,k_2)=(l,e)=n$, similar results to Theorems \ref{thmA} and \ref{thmB} hold for the quotient manifolds $(S^2\times S^3)/\bbz_n$ and $X_\infty/\bbz_n$. 
\end{remark}

\section{Concluding Remarks}

Recall \cite{BG05} that a Sasakian structure $\cals=(\xi,\eta,\Phi,g)$ is said to be {\it positive (negative)} if its basic first Chern class $c_1(\calf_\xi)$ can be represented by a positive (negative) definite $(1,1)$-form. It is {\it null} if $c_1(\calf_\xi)=0$, and {\it indefinite} otherwise. The following result follows directly from Theorem 8.1.14 of \cite{BG05} and Proposition 4.4 of \cite{BGS06}.

\begin{lemma}\label{CRauttype}
Let $\cals=(\xi,\eta,\Phi,g)$ be a Sasakian structure with underlying CR structure $(\cald,J)$. Suppose that $\dim \grk(\cald,J)>1$. Then the type of $\cals$ is either positive or indefinite.
\end{lemma}

In a forthcoming work \cite{BoPa10} it will be seen that generally the type is not an invariant of the Sasaki cone. Here it is easy to see from the constructions above that the type is not an invariant of the bouquet. For example consider the contact structure $\cald_{5,1}$ on $S^2\times S^3$. There are five Sasaki cones $\grk(\cald_{5,1},J_{2m})$ where $m=0,\cdots,4$. We can label the corresponding extremal Sasakian structures as $\cals_{2m}=(\xi,\eta_{5,1},\Phi_{2m},g_{2m})$ where $J_{2m}=\Phi_{2m} |_{\cald_{5,1}}$ and $g_{2m}= d\eta_{5,1}\circ(\Phi_{2m}\otimes \BOne)+\eta_{5,1}\otimes \eta_{5,1}$. The Sasakian structure $\cals_0$ is positive and has constant scalar curvature; whereas, the others $\cals_{2m}$ with $m=1,2,3,4$ are indefinite with non-constant scalar curvature. Similarly, the contact structure $\tcald_{5,4}$ on $X_\infty$ has four Sasaki cones $\grk(\tcald_{4,1},J_{2m+1})$ with $m=0,1,2,3$. The corresponding extremal Sasakian structures are $\cals_{2m+1}=(\xi,\teta_{4,1},\Phi_{2m+1},g_{2m+1})$ where $J_{2m+1}=\Phi_{2m+1} |_{\tcald_{4,1}}$ and $g_{2m+1}= d\teta_{5,1}\circ(\Phi_{2m+1}\otimes \BOne)+\teta_{4,1}\otimes \teta_{4,1}$. Again $\cals_1$ is positive, and the others are indefinite, but now they are all extremal of non-constant scalar curvature.

The case $(l,e)=(2,1)$ and $m=0$ presents an interesting case even though the Sasaki bouquet degenerates to a single Sasaki cone. Here we have $c_1((\tcald_{2,1},J_0))=0$ and $M=S^2\times S^3$. So by Calabi the induced Sasaki structure over the complex manifold $\widetilde{\bbc\bbp}^2$ has an extremal Sasaki metric with non-constant scalar curvature. However, it was shown in \cite{MaSp06} that 
the contact structure $\tcald_{2,1}$ admits an irregular Sasaki-Einstein structure, which, of course, is extremal with constant scalar curvature. In this case it would be interesting to see how big the extremal set $\ge(\tcald_{l,e},J_0)$ is. This is actually a special case of a much more general result of Futaki, Ono, and Wang \cite{FOW06} which says that for any positive toric Sasakian structure with $c_1(\cald)=0$ there is a deformation to a Reeb vector field in the Sasaki cone whose Sasakian metric is Sasaki-Einstein. This discussion begs the questions:

\noindent {\it Does every toric contact structure of Reeb type on $S^2\times S^3$ or $X_\infty$ have a constant scalar curvature Sasakian metric somewhere in its Sasaki cone?, or more generally for any toric contact structure of Reeb type?}

Another observation comes from the well known fact \cite{MoKo06} that for each $m>0$ there is a complex analytic family of even Hirzebruch surfaces $S_{2l}(t)$ satisfying $S_{2l}(t)\approx S_{2l}$ for $t\neq 0$ and $m>l$, but that $S_{2l}(0)\approx S_{2m}.$ This implies that the moduli space of complex structures on $S^2\times S^2$ is non-Hausdorff. A similar result holds for the moduli space of complex structures on $\widetilde{\bbc\bbp}^2$. These results imply that the moduli space of extremal Sasakian structures with underlying contact structures $\cald_{k_1,k_2}$ on $S^2\times S^3$ or $\tcald_{l,e}$ on $X_\infty$ is non-Hausdorff as well.  

Note added: The first question above has been recently answered in the affirmative by \'Eveline Legendre \cite{Leg10}. Moreover, she shows that if $k_1>5k_2$ the toric contact structure $\cald_{k_1,k_2}$ on $S^2\times S^3$ has two distinct rays in the Sasaki cone with constant scalar curvature Sasakian metrics. So unlike the Sasaki-Einstein case \cite{CFO07,FOW06} constant scalar curvature rays are not unique in a Sasaki cone.

\section*{Acknowledgements}
I would like to thank V. Apostolov for many discussions about extremal K\"ahler metrics.

\def\cprime{$'$} \def\cprime{$'$} \def\cprime{$'$} \def\cprime{$'$}
  \def\cprime{$'$} \def\cprime{$'$} \def\cprime{$'$} \def\cprime{$'$}
  \def\cdprime{$''$} \def\cprime{$'$} \def\cprime{$'$}
\providecommand{\bysame}{\leavevmode\hbox to3em{\hrulefill}\thinspace}
\providecommand{\MR}{\relax\ifhmode\unskip\space\fi MR }
% \MRhref is called by the amsart/book/proc definition of \MR.
\providecommand{\MRhref}[2]{%
  \href{http://www.ams.org/mathscinet-getitem?mr=#1}{#2}
}
\providecommand{\href}[2]{#2}

\end{document}